# NONQUADRATIC ESTIMATORS OF A QUADRATIC FUNCTIONAL[1]

By T. Tony Cai and Mark G. Low

*University of Pennsylvania*

Estimation of a quadratic functional over parameter spaces that are not quadratically convex is considered. It is shown, in contrast to the theory for quadratically convex parameter spaces, that optimal quadratic rules are often rate suboptimal. In such cases minimax rate optimal procedures are constructed based on local thresholding. These nonquadratic procedures are sometimes fully efficient even when optimal quadratic rules have slow rates of convergence. Moreover, it is shown that when estimating a quadratic functional nonquadratic procedures may exhibit different elbow phenomena than quadratic procedures.

**1. Introduction.** The Gaussian sequence model

$$(1) \qquad Y_i = \theta_i + n^{-1/2} z_i, \qquad i = 1, 2, \ldots,$$

where $z_i$ are i.i.d. standard normal random variables, often serves as a general prototypical model in nonparametric function estimation settings. For example, it is exactly equivalent to a white noise with drift model and can also be used to approximate nonparametric regression and density estimation models. For the sequence model considerable attention has focused on estimating linear and nonlinear functionals of the infinite dimensional mean vector $\theta = (\theta_1, \theta_2, \ldots)$.

One particularly important nonlinear functional is the quadratic functional $Q(\theta) = \sum_{i=1}^{\infty} \theta_i^2$. Early results on this and related problems were given in [3, 11, 15, 17, 18, 19]. More recent results can be found in [13, 22, 23].

The problem of estimating this quadratic functional is closely connected to the construction of confidence balls in nonparametric function estimation.

Received February 2004; revised December 2004.
[1] Supported in part by NSF Grant DMS-03-06576.

See, for example, [6, 12, 14, 25]. In addition, as shown in [3, 11, 15] this problem connects the nonparametric and semiparametric literatures.

One of the interesting features of the quadratic functional estimation problem is that the usual information bound over any bounded subset of $l_2$, as, for example, given in [2], is strictly positive and finite for $\theta \neq 0$. However this bound may or may not be useful. Bickel and Ritov [3] and Ritov and Bickel [26] showed in the context of i.i.d. data that in some cases the information bound is sharp whereas in other cases the information bound is not informative because the minimax rate of convergence is slower than the usual parametric rate. This phenomenon is often known as the elbow phenomenon.

Donoho and Nussbaum [11] and Fan [15] further developed this theory for orthosymmetric quadratically convex parameter spaces such as hyperrectangles or Sobolev balls. In particular, the minimax theory was fully developed in these cases. The elbow phenomenon also occurs in these more general settings. Moreover, quadratic rules occupy a particularly important position in this theory: simple quadratic rules can always be constructed which are minimax rate optimal.

In this paper we focus on the problem of estimating the quadratic functional $Q(\theta)$ over parameter spaces which are not quadratically convex where, as we shall show, quadratic rules are no longer sufficient for minimax estimation. In particular, we explore when the information bound is sharp and when nonquadratic rules are needed to attain the bound. An estimator is called fully efficient if it attains the information bound asymptotically and we say that fully efficient estimation is possible when such an estimator exists. We also consider specific examples of parameter spaces which are not quadratically convex, namely Besov balls $B_{p,q}^\alpha(M)$ and $L_p$ balls $L_p(\alpha, M)$ with $p < 2$. These spaces, defined in Section 2, provide a rich collection of possible parameter spaces. For these spaces we characterize the elbow phenomenon for the performance of optimal quadratic procedures and that of general minimax procedures. In particular, we show that over these spaces when the optimal quadratic procedure does not attain the usual parametric rate minimax rate optimal rules must be nonquadratic.

The paper is organized as follows. In Section 2 we first consider the performance of quadratic procedures over general orthosymmetric parameter spaces. It is known that when the parameter space is quadratically convex optimal quadratic procedures are near minimax. Such an analysis has however not been given for parameter spaces that are not quadratically convex. In fact, as we show, the near minimaxity of optimal quadratic rules typically does not hold when the parameter space is not quadratically convex. It is shown that the maximum risk of quadratic procedures over any parameter space is equal to the maximum risk over the quadratic convex hull. It also follows from our results that for Besov balls and $L_p$ balls with $p < 2$ quadratic rules can be minimax rate optimal only if the minimax quadratic risk is of order $n^{-1}$. For Besov balls and $L_p$ balls the minimax quadratic risk also exhibits the well-known elbow phenomenon. We show that there is a fully efficient quadratic procedure as long as $\alpha > \frac{1}{p} - \frac{1}{4}$ whereas if $\alpha \leq \frac{1}{p} - \frac{1}{4}$ optimal quadratic rules have maximum risk of order $n^{-8s/(1+4s)}$ where $s = \alpha + \frac{1}{2} - \frac{1}{p}$.

In Section 3 we develop nonquadratic procedures for estimating the quadratic functional $Q(\theta)$ over Besov balls and $L_p$ balls with $p < 2$. We show that optimal nonquadratic procedures exhibit a different elbow phenomenon. A local thresholding estimator is constructed and is shown to be fully efficient over Besov balls and $L_p$ balls with $p < 2$ and $\alpha > \frac{1}{2p}$. Hence when $p < 2$ and $\frac{1}{2p} < \alpha \leq \frac{1}{p} - \frac{1}{4}$ there are fully efficient nonquadratic estimators while all quadratic rules are rate suboptimal.

Section 3 also considers estimating $Q(\theta)$ over Besov balls and $L_p$ balls with $p < 2$ and $\alpha \leq \frac{1}{2p}$. In this case it is shown that the minimax rate of convergence is $n^{-(2-p/(1+2ps))}$ where $s = \alpha + \frac{1}{2} - \frac{1}{p}$, and hence optimal quadratic rules are once again suboptimal since $2 - \frac{p}{1+2ps} > \frac{8s}{1+4s}$. A nonquadratic estimator is constructed which has risk within a constant factor of the minimax risk.

A distinct feature of the case $p < 2$ is that the hardest hyperrectangle submodel is not as difficult as the full model. In contrast, in the dense case of $p \geq 2$ hyperrectangle submodels can be chosen which yield not only useful lower bounds but also lead to rate optimal quadratic procedures. See [11]. For $p < 2$ the worst case can be captured by a mixture prior supported on a large collection of hyperrectangles. Lower bounds are developed in Section 3.3 based on this mixture prior. Local thresholding procedures which capture any large coefficients are shown to be within a constant factor of these lower bounds.

Section 4 briefly considers the adaptation problem for some special cases. Attention is focused only on adaptive estimation across a collection of parameter spaces over which the minimax rates of convergence are equal. In particular, for the collection of all Besov spaces for which fully efficient estimation is possible a procedure based on term by term thresholding is constructed and is shown to be simultaneously fully efficient over every parameter space in this collection. On the other hand, for a fixed nonparametric rate of convergence another estimator is constructed which is simultaneously rate optimal over all Besov spaces with that given minimax rate of convergence. The general case of adaptation over parameter spaces with different minimax rates of convergence is an interesting but challenging problem. A complete treatment is given in [7].

Connections between the problems of estimating quadratic functionals and a corresponding testing problem is made in Section 5. This testing problem was first studied in [20] and Lepski and Spokoiny [24] developed minimax tests for Besov spaces with $\frac{1}{\alpha} < p < 2$. We show that results developed for the estimation problem in Section 3 extend the theory of testing to cases not previously considered.

Section 6 is devoted to a discussion of connections with other related nonparametric function estimation problems, namely those of global estimation under sum of squared error loss and estimating linear functionals. For example, in global estimation it is known that simple thresholding procedures can yield minimax rate optimal procedures over spaces where a few relatively large coefficients may otherwise lead to a large bias. Proofs are given in Section 7.

**2. Performance of quadratic procedures.** As mentioned in the Introduction, quadratic procedures have received particular attention in the theory of estimating quadratic functionals. They have been shown to work well when the parameter space is orthosymmetric and quadratically convex. Most common parameter spaces, such as Besov balls and $L_p$ balls, are orthosymmetric. In particular, unconditional bases such as wavelet bases transform common function spaces into an orthosymmetric sequence space. See [28]. However, many of these spaces are not quadratically convex. In such cases the performance of quadratic rules has not been studied. In this section we study the performance of quadratic procedures over general orthosymmetric parameter spaces. In addition, we consider in detail estimation over Besov balls and $L_p$ balls with $p < 2$.

2.1. *General orthosymmetric parameter spaces.* Before studying the performance of quadratic procedures over general orthosymmetric parameter spaces it is convenient to introduce some notation. Write $\mathcal{Q}$ for the collection of all quadratic rules, namely those of the form

(2) $$\hat{Q} = \sum a_{i,j} Y_i Y_j + c.$$

Also write $\mathcal{Q}_D$ for the subclass of diagonal quadratic rules, namely those of the form $\sum a_i Y_i^2 + c$. A parameter space $\Theta$ is called orthosymmetric if $\theta = (\theta_1, \theta_2, \ldots, \theta_m, \ldots) \in \Theta$ implies that $(\pm\theta_1, \pm\theta_2, \ldots, \pm\theta_m, \ldots) \in \Theta$ for any choices of the signs $\pm$. An orthosymmetric set $\Theta$ is called quadratically convex if the set $\{(\theta_i^2)_{i=1}^\infty : \theta \in \Theta\}$ is convex.

Write the minimax risk for estimating $Q(\theta) = \sum \theta_i^2$ as

(3) $$R^*(n, \Theta) = \inf_{\hat{Q}} \sup_{\theta \in \Theta} E_\theta (\hat{Q} - Q(\theta))^2$$

and the minimax quadratic risk and minimax diagonal quadratic risk as

(4) $$R_Q^*(n, \Theta) = \inf_{\hat{Q} \in \mathcal{Q}} \sup_{\theta \in \Theta} E_\theta (\hat{Q} - Q(\theta))^2 \quad \text{and}$$
$$R_{DQ}^*(n, \Theta) = \inf_{\hat{Q} \in \mathcal{Q}_D} \sup_{\theta \in \Theta} E_\theta (\hat{Q} - Q(\theta))^2.$$

The problem of estimating quadratic functionals has usually assumed that the parameter space is both orthosymmetric and quadratically convex. Orthosymmetry allows a minimax analysis of general quadratic rules to focus on diagonal quadratic rules. More specifically, for any quadratic rule, say $\hat{Q} = \sum a_{i,j} Y_i Y_j + c$, define $\hat{Q}' = \sum a_{i,i} Y_i^2 + c$. Fan [15] showed that for any orthosymmetric set $\Theta$

$$\sup_{\theta \in \Theta} E_\theta(\hat{Q} - Q(\theta))^2 \geq \sup_{\theta \in \Theta} E_\theta(\hat{Q}' - Q(\theta))^2. \tag{5}$$

In particular, it follows that $R_Q^*(n, \Theta) = R_{DQ}^*(n, \Theta)$ when $\Theta$ is orthosymmetric.

For the analysis of quadratic procedures over general orthosymmetric parameter spaces it is convenient and natural to introduce the quadratic convex hull. For an orthosymmetric set $\Theta$, the quadratic convex hull of $\Theta$ is defined as

$$\text{Q.Hull}(\Theta) = \{(\theta_i)_{i=1}^\infty : (\theta_i^2)_{i=1}^\infty \in \text{Hull}(\Theta_+^2)\}, \tag{6}$$

where $\Theta_+^2 = \{(\theta_i^2)_{i=1}^\infty : (\theta_i)_{i=1}^\infty \in \Theta, \theta_i \geq 0 \ \forall i\}$ and $\text{Hull}(\Theta_+^2)$ denotes the convex hull of the set $\Theta_+^2$. The following theorem characterizes the performance of quadratic rules over an orthosymmetric parameter space.

THEOREM 1. *Let $\hat{Q} \in \mathcal{Q}_D$ be a diagonal quadratic estimator of $Q(\theta) = \sum \theta_i^2$. Then for any orthosymmetric $\Theta$,*

$$\sup_{\theta \in \Theta} E_\theta(\hat{Q} - Q(\theta))^2 = \sup_{\theta \in \text{Q.Hull}(\Theta)} E_\theta(\hat{Q} - Q(\theta))^2. \tag{7}$$

*Consequently the minimax quadratic risk over an orthosymmetric set $\Theta$ equals the minimax quadratic risk over the quadratic convex hull of $\Theta$, that is,*

$$R_Q^*(n; \Theta) = R_Q^*(n; \text{Q.Hull}(\Theta)) = R_{DQ}^*(n; \text{Q.Hull}(\Theta)). \tag{8}$$

Theorem 1 shows that the performance of the optimal quadratic procedure is captured by the minimax quadratic risk over the quadratic convex hull of the parameter space $\Theta$. If in addition $\text{Q.Hull}(\Theta)$ is norm bounded in $l_2$ and convex it follows from Donoho and Nussbaum [11] that $R_Q^*(n; \text{Q.Hull}(\Theta)) \asymp R^*(n; \text{Q.Hull}(\Theta))$ and hence $R_Q^*(n; \Theta) \asymp R^*(n; \text{Q.Hull}(\Theta))$. When $\Theta$ is not quadratically convex $\text{Q.Hull}(\Theta)$ is larger than $\Theta$ and in some cases, as we shall discuss below, $R^*(n; \text{Q.Hull}(\Theta)) \gg R^*(n; \Theta)$. Consequently the optimal quadratic procedure can sometimes have a slower rate of convergence than the minimax rate. As we shall show, such is the case for certain Besov balls and $L_p$ balls.

2.2. *Besov balls and $L_p$ balls.* We now consider as an example Besov balls and $L_p$ balls. The $L_p$ balls are defined by

$$L_p(\alpha, M) = \left\{\theta : \left(\sum i^{ps}|\theta_i|^p\right)^{1/p} \leq M\right\}, \tag{9}$$

where $s = \alpha + \frac{1}{2} - \frac{1}{p} > 0$. Besov balls in sequence space are typically defined in terms of a doubly indexed sequence $\{\theta_{j,k} : j = 0, 1, \ldots, k = 0, \ldots, 2^j - 1\}$. The Besov balls are then defined by

$$B_{p,q}^\alpha(M) = \left\{\theta : \left(\sum_{j=0}^\infty \left(2^{js}\left(\sum_{k=0}^{2^j-1}|\theta_{j,k}|^p\right)^{1/p}\right)^q\right)^{1/q} \leq M\right\}, \tag{10}$$

where $s = \alpha + \frac{1}{2} - \frac{1}{p} > 0$. So that we can give a unified treatment of Besov balls and $L_p$ balls it is convenient for Besov balls to set $\theta_i = \theta_{j,k}$ where $i = 2^j + k$. Noisy observation of Besov coefficients can then still be written as in (1). This convention is used throughout the paper. In addition we shall assume throughout the paper that $p, q, \alpha, s > 0$.

Previous literature has focused primarily on quadratically convex parameter spaces such as Besov balls $B_{p,q}^\alpha(M)$ and $L_p$ balls $L_p(\alpha, M)$ with $p \geq 2$. In particular, Fan [15] gave an analysis for $L_p$ balls with $p \geq 2$ which shows that for the parameter space $\Theta = L_p(\alpha, M)$ the minimax risk satisfies

$$\inf_{\hat{Q}} \sup_{\theta \in \Theta} E_\theta(\hat{Q} - Q(\theta))^2 \asymp n^{-r(\alpha)}, \tag{11}$$

where $r(\alpha) = 1$ when $\alpha \geq \frac{1}{4}$ and $r(\alpha) = \frac{8\alpha}{4\alpha+1}$ when $\alpha < \frac{1}{4}$. An entirely analogous analysis yields the same result when $\Theta = B_{p,q}^\alpha(M)$ for $p \geq 2$. Moreover, Fan [15] gave simple quadratic estimators attaining these minimax rates of convergence over $L_p$ balls. Estimating quadratic functionals over Besov spaces was also considered in [23] where the focus was on adaptive estimation of more general quadratic functionals using model selection.

As pointed out in [11] and [15], one important aspect of the quadratically convex $L_p$ balls is that the difficulty of estimating a quadratic functional is then captured by the hardest hyperrectangle subproblem. This reduction is instrumental in developing a sharp lower bound as well as in the construction of the optimal quadratic rule.

Our focus is on Besov balls and $L_p$ balls with $p < 2$, in which case the parameter spaces are no longer quadratically convex. The standard technique of finding the hardest hyperrectangle subproblem is then no longer sufficient. In fact, quadratic rules are in general suboptimal and the hardest hyperrectangle subproblem need not be as difficult as the full model. Nevertheless the performance of optimal quadratic rules is easy to characterize by the results given in Theorem 1 and an understanding of the quadratic convex hulls of general Besov balls and $L_p$ balls. In fact, when $p < 2$ it is easy to check that

$$\text{(12)} \qquad \text{Q.Hull}\,(L_p(\alpha, M)) = L_2\left(\alpha + \frac{1}{2} - \frac{1}{p}, M\right).$$

See [10]. Similarly it is easy to check that

$$\text{(13)} \qquad \text{Q.Hull}(B_{p,q}^\alpha(M)) = B_{2,q}^{\alpha+1/2-1/p}(M).$$

Write $r^*(\Theta)$ for the exponent $a$ whenever $R^*(n, \Theta) \asymp n^{-a}$ and similarly $r_Q^*(\Theta)$ for the exponent $b$ whenever $R_Q^*(n, \Theta) \asymp n^{-b}$. The following result is then a direct consequence of (11)–(13) and Theorem 1.

COROLLARY 1. *Let $0 < p < 2$. Then $r_Q^*(L_p(\alpha, M)) = r_Q^*(B_{p,q}^\alpha(M)) = \min\{1, \frac{8s}{4s+1}\}$, or equivalently,*

$$\text{(14)} \quad r_Q^*(L_p(\alpha, M)) = \begin{cases} r_Q^*(B_{p,q}^\alpha(M)) = 1, & \text{when } \alpha \geq \frac{1}{p} - \frac{1}{4}, \\ r_Q^*(B_{p,q}^\alpha(M)) = \frac{8s}{4s+1}, & \text{when } \alpha < \frac{1}{p} - \frac{1}{4}. \end{cases}$$

The corollary clearly shows the elbow phenomenon for the minimax quadratic rate of convergence. There is a break between the usual parametric rate of convergence and slower rates of convergence at $\alpha = \frac{1}{p} - \frac{1}{4}$. We shall show later that the break for the minimax risk for nonquadratic procedures is at a smaller value of $\alpha$. This is illustrated in Figure 1 for the case of $p = 1.25$.

When $p < 2$ and $\alpha > \frac{1}{p} - \frac{1}{4}$ it is in fact possible to find a simple procedure which is efficient, asymptotically attaining the exact minimax risk. Let $m = \frac{n}{\log n}$ and set

$$\text{(15)} \qquad \hat{Q}_1 = \sum_{i=1}^{m} \left(Y_i^2 - \frac{1}{n}\right).$$

Then simple calculations and lower bounds given in Section 3 yield

$$\text{(16)} \qquad \sup_{\theta \in \Theta} E_\theta(\hat{Q}_1 - Q(\theta))^2 = R^*(n, \Theta)(1 + o(1)) = \frac{4M^2}{n}(1 + o(1)),$$

where $\Theta = B_{p,q}^\alpha(M)$ or $\Theta = L_p(\alpha, M)$ with $p < 2$ and $\alpha > \frac{1}{p} - \frac{1}{4}$.

**3. Nonquadratic estimators.** In this section we focus on the construction of a new class of nonquadratic estimators which significantly outperforms the optimal quadratic rules for Besov balls and $L_p$ balls when $p < 2$ and $\alpha \leq \frac{1}{p} - \frac{1}{4}$. In this case the minimax quadratic risk converges more slowly than the minimax risk and the result shows that quadratic rules are far from optimal.

We shall consider two separate cases. In the first the nonquadratic estimator is fully efficient over Besov balls and $L_p$ balls when $p < 2$ and $\frac{1}{2p} < \alpha < \frac{1}{p} - \frac{1}{4}$, whereas the best quadratic estimator does not even achieve the usual parametric rate. In the second case with $p < 2$ and $\alpha \leq \frac{1}{2p}$ the nonquadratic estimator has risk converging faster than the minimax quadratic risk. We then derive minimax lower bounds in Section 3.3 which show that the risk of this nonquadratic estimator is within a constant factor of the lower bound, and hence the estimator is minimax rate optimal.

3.1. *Fully efficient estimation*: *Besov and $L_p$ balls with $\frac{1}{2p} < \alpha \leq \frac{1}{p} - \frac{1}{4}$.* In parametric problems, Fisher Information provides a standard benchmark for the performance of an estimator. These bounds are often asymptotically attainable. The information bound is often useful in semiparametric models as well. See, for example, [2]. The problem of estimating a quadratic functional received attention by Ritov and Bickel [26] as an example where the information is strictly positive although it is not always possible to achieve the information bound. In the present context of estimating the quadratic functional $Q(\theta)$ the information can easily be calculated to be $I(\theta) = \frac{n}{4 \sum \theta_i^2}$. Standard theory then yields the lower bound

$$(17) \qquad \inf_{\hat{Q}} \sup_{N_\varepsilon(\theta)} E_\theta(\hat{Q} - Q(\theta))^2 \geq \frac{1}{I(\theta)}(1 + o(1)),$$

where $N_\varepsilon(\theta) = \{(1-t)\theta : 0 \leq t \leq \varepsilon\}$ and $0 < \varepsilon < 1$. It then directly follows that (17) provides a lower bound for the minimax risk over a parameter space $\Theta$ whenever $N_\varepsilon(\theta) \subset \Theta$. In particular, the information bound given in (17) immediately yields

$$(18) \qquad \inf_{\hat{Q}} \sup_{\theta \in \Theta} E_\theta(\hat{Q} - Q(\theta))^2 \geq \frac{4M^2}{n}(1 + o(1))$$

for $\Theta = B_{p,q}^\alpha(M)$ or $\Theta = L_p(\alpha, M)$. In Section 2 a simple quadratic procedure was given which attains the bound given in (18) over Besov and $L_p$ balls with $p < 2$ and $\alpha > \frac{1}{p} - \frac{1}{4}$.

We now consider Besov balls and $L_p$ balls where $p < 2$ and $\frac{1}{2p} < \alpha < \frac{1}{p} - \frac{1}{4}$. Corollary 1 shows that in this case the exponent of the minimax quadratic rate of convergence is $\frac{8s}{4s+1} < 1$. We shall show that in this case fully efficient estimation is possible by using nonquadratic rules. One such fully efficient rule can be given as follows.

Let $m$ be a given positive integer. Divide the indices $i$ beyond $m$ into blocks of increasing block size so that the $j$th block is of the size $2^j m$. For

$i$ in block $j$, set $\tau_i = 2j$, that is,

$$\tag{19} \tau_i = 2\left\lceil \log_2 \frac{i}{m} \right\rceil, \qquad i > m,$$

where $\lceil x \rceil$ denotes the smallest integer greater than or equal to $x$.

For $i \geq m+1$, set $\mu_{n,i} = E_0\{(Y_i^2 - \frac{\tau_i}{n})_+\}$ where the expectation is taken under $\theta = 0$. Let $J_*$ be the largest integer such that $2^{J_*} m \leq n^{1/(4s)} \log n$ where once again $s = \alpha + \frac{1}{2} - \frac{1}{p}$. Set the estimator of the quadratic functional $Q = \sum_{i=1}^{\infty} \theta_i^2$ as

$$\tag{20} \hat{Q}(m) = \sum_{i=1}^{m} \left( Y_i^2 - \frac{1}{n} \right) + \sum_{i=m+1}^{2^{J_*} m} \left\{ \left( Y_i^2 - \frac{\tau_i}{n} \right)_+ - \mu_{n,i} \right\}.$$

The parameter $m$ serves as a tuning parameter. We shall choose different $m$ for different cases. The nonquadratic estimator $\hat{Q}(m)$ is built from a quadratic part and coordinate-wise thresholding with slowly growing threshold levels. The thresholding terms are used to guard against individual large terms in the tail.

For the case $p < 2$ and $\frac{1}{2p} < \alpha < \frac{1}{p} - \frac{1}{4}$, set $m_2 = \frac{n}{\log n}$ in (20) and define the estimator $\hat{Q}_2$ as

$$\tag{21} \hat{Q}_2 = \hat{Q}(m_2) = \sum_{i=1}^{m_2} \left( Y_i^2 - \frac{1}{n} \right) + \sum_{i=m_2+1}^{2^{J_*} m_2} \left\{ \left( Y_i^2 - \frac{\tau_i}{n} \right)_+ - \mu_{n,i} \right\}.$$

The following theorem shows that the estimator $\hat{Q}_2$ is fully efficient.

THEOREM 2. *Let $0 < p < 2$ and $\alpha > \frac{1}{2p}$. Then the estimator $\hat{Q}_2$ defined in* (21) *is fully efficient over Besov balls and $L_p$ balls, that is, it satisfies*

$$\tag{22} \sup_{\theta \in \Theta} E_\theta(\hat{Q}_2 - Q(\theta))^2 = \frac{4M^2}{n}(1 + o(1)),$$

*where $\Theta = B_{p,q}^\alpha(M)$ or $\Theta = L_p(\alpha, M)$.*

Comparing (22) with (14) shows that in the case $\frac{1}{2p} < \alpha \leq \frac{1}{p} - \frac{1}{4}$ nonquadratic rules can be fully efficient although all quadratic rules are necessarily rate suboptimal.

REMARK 1. Note that the condition $s = \alpha + \frac{1}{2} - \frac{1}{p} > 0$ implies that $\alpha > \frac{1}{2p}$ whenever $0 < p \leq 1$. Hence fully efficient estimation of $Q(\theta)$ over Besov balls and $L_p$ balls is always possible when $0 < p < 1$. For $L_p$ balls this has already been noted in [23].

REMARK 2. Although the primary focus in the construction of $\hat{Q}_2$ is on the case $\frac{1}{2p} < \alpha \leq \frac{1}{p} - \frac{1}{4}$, the estimator $\hat{Q}_2$ is also fully efficient when $\alpha > \frac{1}{p} - \frac{1}{4}$. The quadratic part of $\hat{Q}_2$ equals $\hat{Q}_1$ given in (15). The contribution of the thresholding part of $\hat{Q}_2$ is negligible in the case $\alpha > \frac{1}{p} - \frac{1}{4}$.

3.2. *Besov balls and $L_p$ balls with $\alpha \leq \frac{1}{2p}$.* So far we have focused on parameter spaces where fully efficient estimation is possible. We now turn to both Besov balls and $L_p$ balls with $\alpha \leq \frac{1}{2p}$ and construct a nonquadratic estimator which has a much faster rate of convergence than the minimax quadratic rate given in Section 2. This result again shows that quadratic rules are rate suboptimal and there is much to be gained by using nonquadratic rules.

Let $m_3 = n^{p(1+2ps)}$ and set the estimator $\hat{Q}_3$ of the quadratic functional $Q = \sum_{i=1}^{\infty} \theta_i^2$ as $\hat{Q}(m)$ in (20) with $m = m_3$. That is,

$$(23) \quad \hat{Q}_3 = \hat{Q}(m_3) = \sum_{i=1}^{m_3} \left(Y_i^2 - \frac{1}{n}\right) + \sum_{i=m_3+1}^{2^{J_*}m_3} \left\{ \left(Y_i^2 - \frac{\tau_i}{n}\right)_+ - \mu_{n,i} \right\},$$

where once again $J_*$ is the largest integer such that $2^{J_*}m \leq n^{1/(4s)} \log n$. The following provides an upper bound for the risk of the estimator $\hat{Q}_3$.

THEOREM 3. *Let $0 < p < 2$ and $\alpha \leq \frac{1}{2p}$. The estimator $\hat{Q}_3$ given in (23) satisfies*

$$(24) \quad \sup_{\theta \in \Theta} E_\theta(\hat{Q}_3 - Q(\theta))^2 \leq C n^{-(2-p/(1+2ps))}(1+o(1)),$$

*where $C > 0$ is a constant and $\Theta = B_{p,q}^\alpha(M)$ or $\Theta = L_p(\alpha, M)$.*

It is easy to check that if $p < 2$ then $2 - \frac{p}{1+2ps} > \frac{8s}{4s+1}$. Hence quadratic rules are necessarily rate suboptimal when $p < 2$ and $\alpha \leq \frac{1}{2p}$. In the next section it is shown that no estimator has maximum risk converging faster than $n^{-(2-p/(1+2ps))}$ and thus the estimator $\hat{Q}_3$ is minimax rate optimal.

The analysis of the estimators $\hat{Q}_2$ and $\hat{Q}_3$ relies on a detailed analysis of bias and variance of thresholding estimators for each coordinate. The following lemma may also be of independent interest.

LEMMA 1. *Let $X \sim N(\theta, \frac{1}{n})$ and $\tau \geq 1$. Set $\mu_0 = E_0\{(X^2 - \frac{\tau}{n})_+\}$ where the expectation is taken under $\theta = 0$. Let $\hat{Q} = (X^2 - \frac{\tau}{n})_+ - \mu_0$. Then*

$$(25) \quad |\mu_0| \leq \frac{4}{\sqrt{2\pi}n\tau^{1/2}e^{\tau/2}},$$

$$(26) \quad |E_\theta\hat{Q} - \theta^2| \leq \min\left(\frac{2\tau}{n}, \theta^2\right)$$

and the variance of $\hat{Q}$ satisfies

$$\text{Var}(\hat{Q}) \leq \frac{6\theta^2}{n} + \frac{4\tau^{1/2} + 18}{n^2 e^{\tau/2}}. \tag{27}$$

Combining the results given in Section 2 as well as this section for both quadratic and nonquadratic rules, we can compare the optimal rates of convergence over Besov and $L_p$ balls. Figure 1 gives a comparison for the case of $p = 1.25$ as a function of $\alpha$. It illustrates the different elbow phenomena for the minimax rate of convergence and the minimax quadratic rate of convergence.

3.3. *Minimax lower bounds.* As shown earlier, the information bound given in (18) is sharp for estimating $Q(\theta)$ over Besov balls and $L_p$ balls when $0 < p < 2$ and $\alpha > \frac{1}{2p}$, although sometimes nonquadratic rules are needed to attain the bound. When $0 < p < 2$ and $\alpha < \frac{1}{2p}$ the information bound is no longer attainable. In this section we provide an improved lower bound which shows that the minimax rate of convergence is slower than the usual parametric rate. Furthermore these lower bounds show that the nonquadratic estimator $\hat{Q}_3$ given in (23) is minimax rate optimal.

The derivation of the lower bound given in this section differs from the standard technique of inscribing a hardest hyperrectangle and using the Bayes risk for a prior supported on the hyperrectangle as a lower bound to the minimax risk. The hardest hyperrectangle techniques works when the parameter space is quadratically convex. See, for example, [11] and [15]. However, this technique does not work in our context where the hardest hyperrectangle submodel is not as difficult as the full model. The lower bound given below is based on a mixture prior which mixes over a rich collection of hyperrectangles. The mixing increases the difficulty of the Bayes estimation problem and results in a sharper lower bound.

THEOREM 4. *The minimax risks for estimating the quadratic functional $Q(\theta) = \sum \theta_i^2$ over the Besov balls $B_{p,q}^\alpha(M)$ and $L_p$ balls $L_p(\alpha, M)$ satisfy, for some constant $C > 0$,*

$$\inf_\delta \sup_{\theta \in \Theta} E_\theta(\delta - Q(\theta))^2$$

$$\geq \begin{cases} \frac{4M^2}{n}(1 + o(1)), & \text{when } 0 < p < 2 \text{ and } \alpha > \frac{1}{2p}, \\ Cn^{-(2-p/(1+2ps))}, & \text{when } 0 < p < 2 \text{ and } \alpha \leq \frac{1}{2p}, \end{cases} \tag{28}$$

*where $\Theta = B_{p,q}^\alpha(M)$ or $\Theta = L_p(\alpha, M)$.*

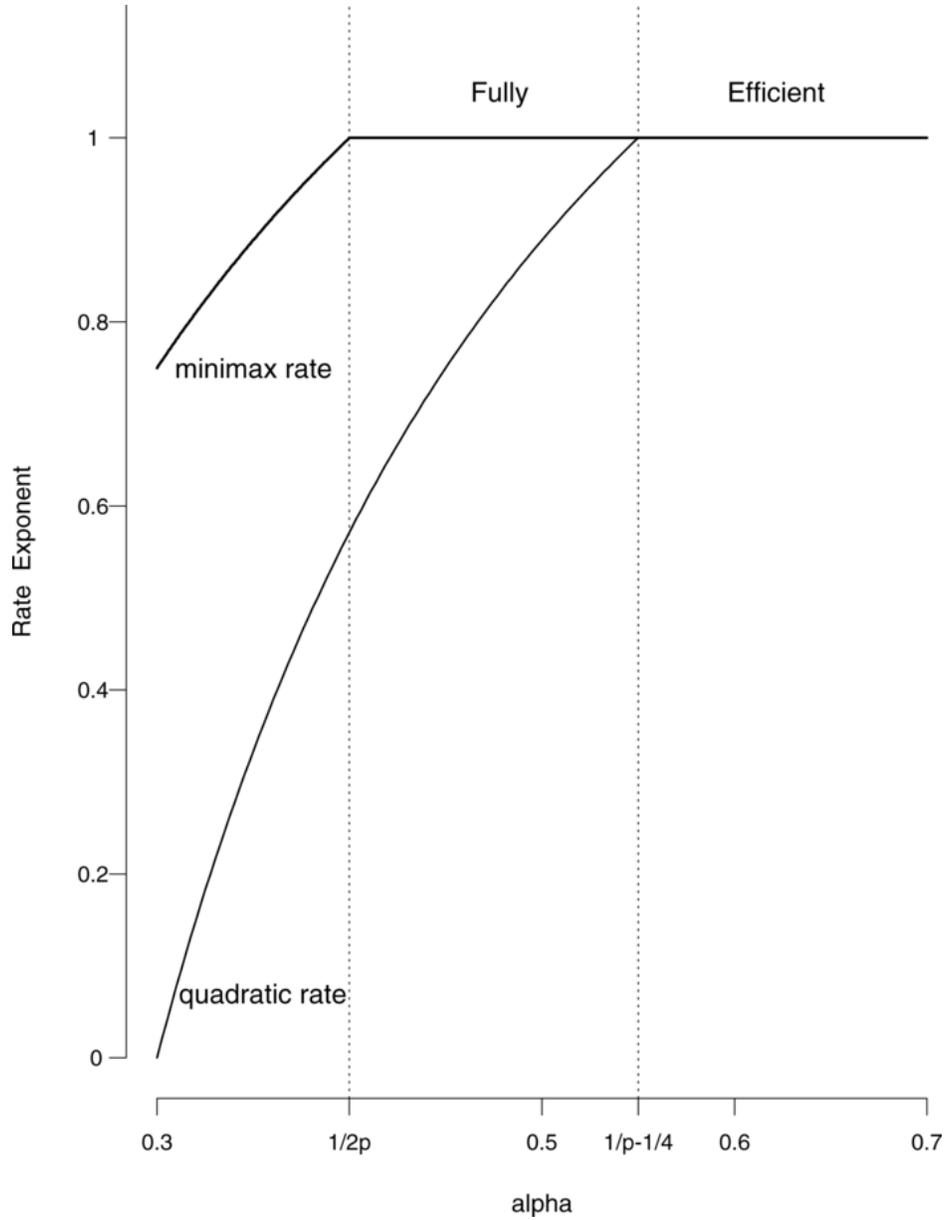

Fig. 1. *Comparison of exponents in the minimax rate of convergence and minimax quadratic rate of convergence for $p = 1.25$.*

The lower bounds show that when $p < 2$ and $\alpha < \frac{1}{2p}$ the optimal rate is slower than the parametric rate. A comparison of the lower bound given above with the upper bound given in (24) shows that in this case the mini-max rate of convergence is $n^{-(2-p/(1+2ps))}$ and the nonquadratic procedure $\hat{Q}_3$ is minimax rate optimal.

TABLE 1
*Comparison of minimax and minimax quadratic rates of convergence*

|  | $0 < p < 2$ | | | $p \geq 2$ | |
| --- | --- | --- | --- | --- | --- |
|  | $\alpha \leq \frac{1}{2p}$ | $\frac{1}{2p} < \alpha \leq \frac{1}{p} - \frac{1}{4}$ | $\alpha > \frac{1}{p} - \frac{1}{4}$ | $\alpha \leq \frac{1}{4}$ | $\alpha > \frac{1}{4}$ |
| $r_Q^*$ | $\frac{8s}{1+4s}$ | $\frac{8s}{1+4s}$ | 1 | $\frac{8\alpha}{1+4\alpha}$ | 1 |
| $r^*$ | $\frac{4ps+2-p}{1+2ps}$ | 1 | 1 | $\frac{8\alpha}{1+4\alpha}$ | 1 |

The results given in Section 2 and this section can be summarized in Table 1 above. For comparison and completeness we also include the well-known results for $p \geq 2$ in the table. As in Section 2 define $r^*$ and $r_Q^*$ to be the exponents of the minimax and minimax quadratic rate of convergence, respectively. We can compare the values of $r^*$ and $r_Q^*$ for all cases in Table 1, where as usual we assume $s = \alpha + \frac{1}{2} - \frac{1}{p} > 0$.

**4. Simple adaptation.** The main focus of the paper is on deficiencies of quadratic estimators and on the minimax performance of the nonquadratic estimators. The estimators $\hat{Q}_2$ and $\hat{Q}_3$ depend on the parameters $\alpha$ and $p$ and are thus not adaptive. A modification of the estimator $\hat{Q}_2$ can achieve full adaptation over collections of Besov balls and $L_p$ balls when fully efficient estimation is possible.

Set $m = \frac{n}{\log n}$ and let $\gamma > 1$ be a constant. Let $J^*$ be the largest integer such that $2^{J^*} m \leq n^\gamma \log n$. Set

$$(29) \qquad \hat{Q}_4 = \sum_{i=1}^{m} \left( Y_i^2 - \frac{1}{n} \right) + \sum_{i=m+1}^{2^{J^*} m} \left\{ \left( Y_i^2 - \frac{\tau_i}{n} \right)_+ - \mu_{n,i} \right\},$$

where $\mu_{n,i}$ is defined the same as in $\hat{Q}_2$. It is then not difficult to show that for all $0 < p < 2$ and $\alpha > \frac{1}{2p} + (\frac{1}{2p} - \frac{1}{2} + \frac{1}{4\gamma})_+$

$$(30) \qquad \sup_{\theta \in \Theta} E_\theta (\hat{Q}_4 - Q(\theta))^2 = \frac{4M^2}{n}(1 + o(1)),$$

where $\Theta = B_{p,q}^\alpha(M)$ or $\Theta = L_p(\alpha, M)$.

Hence, the estimator $\hat{Q}_4$ is adaptively fully efficient over the collection

$$\left\{ B_{p,q}^\alpha(M) : 0 < p < 2, \alpha > \frac{1}{2p} + \left( \frac{1}{2p} - \frac{1}{2} + \frac{1}{4\gamma} \right)_+ \right\}.$$

More interestingly, if we take $J^* = \infty$ and set

$$\hat{Q}_5 = \sum_{i=1}^m \left( Y_i^2 - \frac{1}{n} \right) + \sum_{i=m+1}^\infty \left\{ \left( Y_i^2 - \frac{\tau_i}{n} \right)_+ - \mu_{n,i} \right\}, \tag{31}$$

then the estimator $\hat{Q}_5$ is adaptively fully efficient over all Besov balls and $L_p$ balls with $p < 2$ and $\alpha \geq \frac{1}{2p}$. In fact, $\hat{Q}_5$ is also adaptively fully efficient over Besov balls and $L_p$ balls with $p \geq 2$ and $\alpha > \frac{1}{4}$. It is easy to see from Table 1 that these are the maximum collections of Besov balls and $L_p$ balls over which fully efficient estimation is possible. Therefore the estimator $\hat{Q}_5$ is adaptively fully efficient whenever fully efficient estimation is possible. We summarize the results in the following theorem.

THEOREM 5. *The estimator $\hat{Q}_5$ defined in* (31) *satisfies*

$$\sup_{\theta \in \Theta} E_\theta(\hat{Q}_5 - Q(\theta))^2 = \frac{4M^2}{n}(1 + o(1)), \tag{32}$$

*where $\Theta = B_{p,q}^\alpha(M)$ or $\Theta = L_p(\alpha, M)$, with $0 < p < 2$ and $\alpha > \frac{1}{2p}$ or $p \geq 2$ and $\alpha > \frac{1}{4}$.*

Similarly we can also consider the case where the minimax rate of convergence is nonparametric. Fix a constant $0 < r < 1$ and let

$$\Omega(r) = \left\{ (\alpha, p) : 0 < p < 2, \ 0 < \alpha < \frac{1}{2p}, \ \frac{p}{1 + 2ps} = 2 - r \right\},$$

where, as usual, $s = \alpha + \frac{1}{2} - \frac{1}{p} > 0$. Note that the minimax rate of convergence for estimating the quadratic functional $Q(\theta) = \sum \theta_i^2$ over any Besov or $L_p$ ball with parameters $(\alpha, p) \in \Omega(r)$ is $n^{-r}$.

Let $m = n^{2-r}$ and let $\tau_i$ be defined as in (19). Set

$$\hat{Q}_6 = \sum_{i=1}^m \left( Y_i^2 - \frac{1}{n} \right) + \sum_{i=m+1}^\infty \left\{ \left( Y_i^2 - \frac{\tau_i}{n} \right)_+ - \mu_{n,i} \right\}. \tag{33}$$

Then it is easy to show that the estimator $\hat{Q}_6$ adaptively attains the minimax rate of convergence $n^{-r}$ over each Besov ball $B_{p,q}^\alpha(M)$ or $L_p$ ball $L_p(\alpha, M)$ with $(\alpha, p) \in \Omega(r)$.

The discussion given above is restricted to cases where the minimax rates of convergence over all the parameter spaces in the collection are the same. A more general approach should consider adaptation over spaces with different minimax rates of convergence. The more standard case where $p > 2$ has been considered by Klemelä [22]. The general case is an interesting but challenging problem. A complete treatment is given in [7].

**5. Connection between estimation and testing.** As is common in statistical inference there are strong connections among the problems of estimation and testing quadratic functionals. In this context the testing problem which has received most attention is that of testing the null hypothesis

$$H_0 : \theta = 0$$

against the alternative

$$H_a : \sum \theta_i^2 \geq a_n.$$

The difficulty of this testing problem depends on assumptions imposed on the unknown $\theta$. Particularly important early work on this problem can be traced to Ingster [20].

There are two major related goals in these problems. One goal is to find the test which, given a particular choice of $a_n$, minimizes the sum of the type I and maximum type II errors. Alternatively we may fix the maximal sum of the type I and type II errors and try to find the smallest possible $a_n$ compatible with this constraint along with the corresponding test. More specifically, for a given $0 < \gamma < 1$ let $a_n(\gamma)$ be the smallest choice of $a_n$ for which there is a test with type I plus maximal type II error less than or equal to $\gamma$.

The solution to this testing problem always yields lower bounds to the corresponding estimation problem as follows. First note that every estimator $\hat{Q}$ of $\sum \theta_i^2$ gives rise to a test of this hypothesis in the following way. If $\hat{Q} \leq \frac{a_n(\gamma)}{2}$, then the null hypothesis $H_0$ is accepted and if $\hat{Q} > \frac{a_n(\gamma)}{2}$ the null is rejected. It then immediately follows that

$$(34) \qquad \sup E\left(\hat{Q} - \sum \theta_i^2\right)^2 \geq \frac{1}{2}\gamma \frac{a_n^2(\gamma)}{4} = \frac{1}{8}\gamma a_n^2(\gamma).$$

It is then easy to connect an asymptotic statement about the testing problem into asymptotic lower bounds for the estimation problem. For example, if $r_t$ is the optimal minimax rate for the testing problem, namely $a_n(\gamma) \sim n^{-r_t}$, it then follows that

$$(35) \qquad \inf \sup E\left(\hat{Q} - \sum \theta_i^2\right)^2 \geq C n^{-2r_t}(1 + o(1)).$$

Hence knowledge about the optimal rate in testing immediately yields a lower bound for the optimal rate in the estimating problem. Likewise upper bounds on the estimation problem yield upper bounds on the testing problem. For example, if

$$\sup_{\theta \in \Theta} E\left(\hat{Q} - \sum \theta_i^2\right)^2 \sim n^{-r_e},$$

then

$$a_n(\gamma) \leq C n^{-r_e/2}(1 + o(1)).$$

This testing problem has been considered in [24] over Besov balls $B_{p,q}^\alpha$ for the cases where $\frac{1}{\alpha} < p < 2$. Although the testing theory has not been developed for the cases where $\alpha \leq \frac{1}{p}$, the present paper does immediately yield the optimal rates for testing whenever $\alpha \leq \frac{1}{2p}$. In this case the optimal rate for the testing is $n^{-(1-p/(2(1+2ps)))}$ and this rate has the same functional form as the minimax rate given in [24] in the case $\alpha > \frac{1}{p}$. For the range $\frac{1}{2p} < \alpha \leq \frac{1}{p}$ the estimation problem appears to be "harder" than the testing problem and our results on estimating the quadratic functional do not yield sharp lower bounds or upper bounds on testing.

**6. Discussion.** There are strong similarities between estimating a quadratic functional over a parameter space that is not quadratically convex and that of estimating linear functionals over a nonconvex parameter space. For estimating the quadratic functional it is shown in Section 2 that the maximum risk of quadratic procedures over a parameter space is equal to the maximum risk over the quadratic convex hull of the parameter space. On the other hand it was shown in [5] that for estimating linear functionals the maximum risk of linear procedures over a parameter space is equal to the maximum risk over the convex hull of the parameter space.

There is also some similarity between the work on estimating a quadratic functional and that of estimating all the coefficients under sum of squared error loss. In both problems extra care must be taken for parameter spaces where a few large coefficients can degrade the performance of naive estimators. Under sum of squared error loss the naive estimators correspond to linear estimators. Such estimators can perform well for $p \geq 2$. In particular there are simple linear procedures which are minimax rate optimal. On the other hand, if $p < 2$ then minimax rate optimal procedures must be nonlinear. The case where $p \geq 2$ is sometimes referred to as the dense case since in this situation the difficulty of the problem is caused by situations where there are many small coefficients. The case where $p < 2$ corresponds to sparse situations where there may be a "few" large coefficients which if not estimated well inflate the risk. For global estimation under sum of squared error loss Donoho and Johnstone [8] have shown that fairly simple term-by-term thresholding rules can then yield minimax rate optimal procedures. In density estimation problems Donoho, Johnstone, Kerkyacharian and Picard [9] showed a similar phenomenon exists.

For estimating quadratic functionals, the naive estimators are quadratic rather than linear. Similar to the problem of global estimation the case $p \geq 2$ is easiest. Minimax rate optimal quadratic estimators always exist. However, for estimating a quadratic functional quadratic rules can sometimes be asymptotically fully efficient even in cases when $p < 2$. Such is the case when $p < 2$ and $\alpha > \frac{1}{p} - \frac{1}{4}$. When $p < 2$ and $\alpha < \frac{1}{p} - \frac{1}{4}$, quadratic rules are rate suboptimal while in some cases nonquadratic rules can be fully efficient. As in the case of global estimation term-by-term thresholding can yield minimax rate optimal procedures.

In global estimation minimax rate optimal procedures can be based either on soft thresholding or hard thresholding. The same holds when estimating the quadratic functional. The form of this thresholding is not important. There is an analog of Lemma 1 for hard thresholding and so minimax rate optimal procedures can be based on hard thresholding. More specifically, let the estimator $\tilde{Q}(m)$ be defined as

$$(36) \qquad \tilde{Q}(m) = \sum_{i=1}^{m}\left(Y_i^2 - \frac{1}{n}\right) + \sum_{i=m+1}^{2^{J_*}m}\left\{Y_i^2 I\left(Y_i^2 > \frac{\tau_i}{n}\right) - \rho_{n,i}\right\},$$

where $\tau_i$ is given as in (19) and $\rho_{n,i} = E_0\{Y_i^2 I(Y_i^2 > \frac{\tau_i}{n})\}$ with the expectation taken under $\theta = 0$. Then the results of Theorems 2, 3 and 5 hold for $\tilde{Q}(m)$ with the same choices of $m$.

Finally we should also note that the term-by-term thresholding procedures used here are quite different from the global thresholding rules used in [21] and [27], which are designed for estimation over quadratically convex parameter spaces where the worst case is always given by a large number of small coefficients but where the exact locations of these coefficients are unknown.

**7. Proofs.** For proofs involving both Besov balls and $L_p$ balls we shall only give details for the $L_p$ balls since the proofs for Besov balls are entirely analogous. The proofs of Theorems 2, 3 and 5 all rely on the technical result given in Lemma 1 and are similar. We shall present a detailed proof for Theorem 3, a brief proof for Theorem 2 and omit the proof for Theorem 5. In the proofs we shall denote by $C$ a positive constant not depending on $n$ that may vary from place to place.

PROOF OF THEOREM 1. Since $\Theta \subseteq \text{Q.Hull}(\Theta)$, it is obvious that

$$\sup_{\theta \in \Theta} E_\theta(\hat{Q} - Q(\theta))^2 \le \sup_{\theta \in \text{Q.Hull}(\Theta)} E_\theta(\hat{Q} - Q(\theta))^2.$$

Let $\hat{Q}$ be a diagonal quadratic estimator of $Q(\theta)$. Write $\hat{Q} = \sum_i a_i Y_i^2 + b$. Let $\theta \in \text{Q.Hull}(\Theta)$ and $\theta^2 = \sum_j \lambda_j (\theta^{(j)})^2$ with $\theta^{(j)} \in \Theta$, $\lambda_j \ge 0$ and $\sum_j \lambda_j = 1$. Write

$$V(\theta) = \text{Var}_\theta(\hat{Q}) \quad \text{and} \quad B(\theta) = E_\theta \hat{Q} - Q(\theta).$$

Then

$$V(\theta) = \frac{4\sum_i a_i^2 \theta_i^2}{n} + \frac{2\sum_i a_i^2}{n^2} = \sum_j \lambda_j \left\{\frac{4\sum_i a_i^2 (\theta_i^{(j)})^2}{n} + \frac{2\sum_i a_i^2}{n^2}\right\}$$

$$= \sum_j \lambda_j V(\theta^{(j)})$$

and

$$B(\theta) = \sum_i a_i \theta_i^2 + \sum_i \frac{a_i}{n} + b - \sum_i \theta_i^2$$

$$= \sum_j \lambda_j \left\{ \sum_i a_i (\theta_i^{(j)})^2 + \sum_i \frac{a_i}{n} + b - \sum_i (\theta_i^{(j)})^2 \right\}$$

$$= \sum_j \lambda_j B(\theta^{(j)}).$$

Let $\max_j \{V(\theta^{(j)}) + B^2(\theta^{(j)})\} = D$. Then the Cauchy–Schwarz inequality yields

$$E_\theta(\hat{Q} - Q(\theta))^2 = V(\theta) + B^2(\theta) = \sum_j \lambda_j V(\theta^{(j)}) + \left( \sum_j \lambda_j B(\theta^{(j)}) \right)^2$$

$$\leq D - \sum_j \lambda_j B^2(\theta^{(j)}) + \left( \sum_j \lambda_j B(\theta^{(j)}) \right)^2$$

$$\leq D$$

$$\leq \sup_{\theta' \in \Theta} E_{\theta'}(\hat{Q} - Q(\theta'))^2.$$

Hence for any diagonal quadratic estimator $\hat{Q}$

$$(37) \quad \sup_{\theta \in \mathrm{Q.Hull}(\Theta)} E_\theta(\hat{Q} - Q(\theta))^2 = \sup_{\theta \in \Theta} E_\theta(\hat{Q} - Q(\theta))^2.$$

Since $\Theta$ is orthosymmetric it follows from [15] that minimax quadratic procedures are found within the class of diagonal quadratic procedures. So it follows from (37) that for any quadratic estimator $\hat{Q}$

$$\sup_{\theta \in \mathrm{Q.Hull}(\Theta)} E_\theta(\hat{Q} - Q(\theta))^2 = \sup_{\theta \in \Theta} E_\theta(\hat{Q} - Q(\theta))^2,$$

and hence (8) holds. $\square$

PROOF OF LEMMA 1. Denote by $\phi(z)$ and $\Phi(z)$ the density and cumulative distribution function of $Z$, respectively, and set $\tilde{\Phi}(z) = 1 - \Phi(z)$. It then follows from the alternating series bound for Gaussian tails $\tilde{\Phi}(z) \geq (\frac{1}{z} - \frac{1}{z^3})\phi(z)$ for $z > 0$ that

$$\mu_0 = \frac{2}{n} \int_{\tau^{1/2}}^\infty (z^2 - \tau) \frac{1}{\sqrt{2\pi}} e^{-z^2/2} \, dz$$

$$= \frac{2\tau^{1/2}}{\sqrt{2\pi n}e^{\tau/2}} - \frac{2(\tau-1)}{n}\tilde{\Phi}(\tau^{1/2})$$

(38)
$$\leq \frac{2\tau^{1/2}}{\sqrt{2\pi n}e^{\tau/2}} - \frac{2(\tau-1)}{\sqrt{2\pi n}e^{\tau/2}}\left(\frac{1}{\tau^{1/2}} - \frac{1}{\tau^{3/2}}\right)$$

$$\leq \frac{4}{\sqrt{2\pi n}\tau^{1/2}e^{\tau/2}}.$$

Set $B(\theta) = E_\theta \hat{Q} - \theta^2 = E_\theta(X^2 - \frac{\tau}{n})_+ - \mu_0 - \theta^2$. It is easy to check that

(39) $$\theta^2 - \frac{\tau}{n} \leq E\left(X^2 - \frac{\tau}{n}\right)_+ \leq \theta^2 + \frac{1}{n}.$$

Hence

(40) $$|B(\theta)| \leq \frac{\tau}{n} + \mu_0 \leq \frac{2\tau}{n}.$$

Straightforward calculation yields for $\theta \geq 0$

(41)
$$B'(\theta) = \frac{2}{\sqrt{n}}[\phi(\tau^{1/2} - n^{1/2}\theta) - \phi(\tau^{1/2} + n^{1/2}\theta)]$$
$$- 2\theta[\Phi(\tau^{1/2} - n^{1/2}\theta) - \Phi(-\tau^{1/2} - n^{1/2}\theta)]$$

and

(42)
$$B''(\theta) = 2\tau^{1/2}[\phi(\tau^{1/2} - n^{1/2}\theta) + \phi(\tau^{1/2} + n^{1/2}\theta)]$$
$$- 2[\Phi(\tau^{1/2} - n^{1/2}\theta) - \Phi(-\tau^{1/2} - n^{1/2}\theta)].$$

It suffices to only consider $\theta \geq 0$ since $B(\theta) = B(-\theta)$. It follows immediately from (41) that for all $\theta \geq 0$, $B'(\theta) \geq -2\theta$ and hence

(43) $$B(\theta) \geq -\theta^2.$$

On the other hand, for $0 \leq \theta \leq \frac{1}{\sqrt{n}}$, equation (42) yields

(44) $$B''(\theta) \leq \sup_{\tau \geq 1}\{2\tau^{1/2}[\phi(\tau^{1/2} - 1) + \phi(\tau^{1/2})]\} \leq 2.$$

Note that $B'(0) = 0$ and hence for $0 \leq \theta \leq \frac{1}{\sqrt{n}}$

(45) $$B'(\theta) \leq 2\theta.$$

For $\theta \geq \frac{1}{\sqrt{n}}$ it follows from (41) that

(46) $$B'(\theta) \leq \frac{2}{\sqrt{n}} \leq 2\theta$$

and it follows from $B(0) = 0$ that for all $\theta$, $B(\theta) \leq \theta^2$. Hence for all $\theta$

$$|B(\theta)| \leq \theta^2 \tag{47}$$

and (26) now follows from (40) and (47). The proof will be complete once we establish (27). First we state and prove the following simple lemma.

LEMMA 2. *For any two random variables $X$ and $Y$,*

$$\mathrm{Var}(\max\{X,Y\}) \leq \mathrm{Var}\, X + \mathrm{Var}\, Y. \tag{48}$$

*In particular, for any random variable $X$,*

$$\mathrm{Var}\,((X)_+) \leq \mathrm{Var}\, X. \tag{49}$$

PROOF. Without loss of generality we can assume $\mu_X = 0$ and $\mu_Y \geq 0$. Let $Z = \max\{X, Y\}$. Then

$$EZ^2 \leq EX^2 + EY^2 \tag{50}$$

and

$$EZ \geq \mu_Y. \tag{51}$$

Hence

$$\mathrm{Var}\, Z = EZ^2 - (EZ)^2 \leq EX^2 + EY^2 - \mu_Y^2 = \mathrm{Var}\, X + \mathrm{Var}\, Y. \tag{52} \qquad \square$$

We now turn to the proof of (27). When $\theta^2 \geq \frac{1}{n}$, it follows from Lemma 2 that

$$\mathrm{Var}(\hat{Q}) \leq \mathrm{Var}(X^2) = \frac{4\theta^2}{n} + \frac{2}{n^2} \leq \frac{6\theta^2}{n}$$

and so (27) holds.

Now consider $\theta^2 < \frac{1}{n}$. Because of the symmetry, it suffices to consider $0 \leq \theta < \frac{1}{\sqrt{n}}$. Note that direct calculations show for $0 \leq \theta < \frac{1}{\sqrt{n}}$

$$\mathrm{Var}(\hat{Q}) \leq E\left\{\left(X^2 - \frac{\tau}{n}\right)_+\right\}^2 = \frac{\sqrt{n}}{\sqrt{2\pi}} \int_{x^2 \geq \tau/n} \left(x^2 - \frac{\tau}{n}\right)^2 e^{-n/2(x-\theta)^2}\, dx$$

$$= \frac{1}{\sqrt{2\pi}n^2} \int_{(z+n^{1/2}\theta)^2 \geq \tau} ((z+n^{1/2}\theta)^2 - \tau)^2 e^{-z^2/2}\, dz$$

$$\leq \frac{1}{\sqrt{2\pi}n^2} \int_{\tau^{1/2}-n^{1/2}\theta}^{\tau^{1/2}} ((z+n^{1/2}\theta)^2 - \tau)^2 e^{-z^2/2}\, dz$$

$$+ \frac{2}{\sqrt{2\pi}n^2} \int_{\tau^{1/2}}^{\infty} ((z+1)^2 - \tau)^2 e^{-z^2/2}\, dz$$

NONQUADRATIC ESTIMATORS OF A QUADRATIC FUNCTIONAL 21

$$\leq \frac{1}{\sqrt{2\pi}n^2} \int_{\tau^{1/2}-n^{1/2}\theta}^{\tau^{1/2}} ze^{-z^2/2}\,dz \cdot \sup_{z\in[\tau^{1/2}-n^{1/2}\theta,\tau^{1/2}]} z^{-1}((z+n^{1/2}\theta)^2-\tau)^2$$

$$+ \frac{2}{\sqrt{2\pi}n^2}\int_{\tau^{1/2}}^{\infty}[z^4+4z^3+(-2\tau+6)z^2+(-4\tau+4)z+(\tau-1)^2]$$

$$\times e^{-z^2/2}\,dz$$

$$\equiv V_1 + V_2.$$

Note that

$$\int_{\tau^{1/2}-n^{1/2}\theta}^{\tau^{1/2}} ze^{-z^2/2}\,dz \leq \int_{\tau^{1/2}-1}^{\tau^{1/2}} ze^{-z^2/2}\,dz = e^{-1/2(\tau^{1/2}-1)^2} - e^{\tau/2}$$

and

$$\sup_{z\in[\tau^{1/2}-n^{1/2}\theta,\tau^{1/2}]} z^{-1}((z+n^{1/2}\theta)^2 - \tau)^2 = \sup_{x\in[0,n^{1/2}\theta]} \frac{x(x+2\tau^{1/2})^2}{1+(\tau^{1/2}-n^{1/2}\theta)x^{-1}}$$

$$= n\theta^2 \tau^{-1/2}(2\tau^{1/2}+n^{1/2}\theta)^2$$

$$\leq n\theta^2 \tau^{-1/2}(2\tau^{1/2}+1)^2.$$

Hence

$$V_1 \leq \frac{\theta^2}{n} \cdot \sup_{x\geq 1}\left\{\frac{1}{\sqrt{2\pi}}x^{-1}(2x+1)^2[e^{-1/2(x-1)^2} - e^{x^2/2}]\right\}$$

$$\leq \frac{3\theta^2}{n}.$$

We now turn to the term $V_2$. Note that for any $p \geq 0$

$$\int_a^{\infty} z^p e^{-z^2/2}\,dz = a^{p-1}e^{-a^2/2} + (p-1)\int_a^{\infty} z^{p-2}e^{-z^2/2}\,dz,$$

and in particular for any $a > 0$

$$\int_a^{\infty} e^{-z^2/2}\,dz \leq a^{-1}e^{-a^2/2}.$$

Hence, after some algebra, we have

$$V_2 \leq \frac{10\tau^{1/2}+20\tau^{-1/2}+24}{\sqrt{2\pi}n^2 e^{\tau/2}} \leq \frac{4\tau^{1/2}+18}{n^2 e^{\tau/2}}.$$

Hence, for $0 \leq \theta \leq \frac{1}{\sqrt{n}}$

$$\mathrm{Var}(\hat{Q}) \leq \frac{3\theta^2}{n} + \frac{4\tau^{1/2}+18}{n^2 e^{\tau/2}}$$

and consequently, for all $\theta$, $\mathrm{Var}(\hat{Q}) \leq \frac{6\theta^2}{n} + \frac{4\tau^{1/2}+18}{n^2 e^{\tau/2}}$. □

PROOF OF THEOREM 3. Set $m = m_3 = n^{p/(1+2ps)}$. Note that for $X \sim N(\mu, \sigma^2)$, $\text{Var}(X^2) = 4\mu^2\sigma^2 + 2\sigma^4$. Note also that $\tau_i = 2j$ for all coordinates in the $j$th block beyond the initial $m$ terms. It then follows from Lemma 1 that

$$\text{Var}(\hat{Q}_3) \leq \frac{2m}{n^2} + \frac{4\sum_{i=1}^{m}\theta_i^2}{n} + \frac{6\sum_{i=m+1}^{2^{J_*}m}\theta_i^2}{n}$$

$$(53) \qquad + \sum_{j=1}^{J_*} 2^{j-1}m \cdot \frac{4(2j)^{1/2} + 18}{n^2 e^j}$$

$$\leq Cn^{-(2-p/(1+2ps))}(1+o(1))$$

for some constant $C > 0$, where the last step follows from the fact that for any $b > 0$

$$(54) \qquad \sum_{j=1}^{\infty} j^{1/2} e^{-bj} < \infty.$$

For the bias note that equation (26) in Lemma 1 yields

$$(55) \qquad |\text{Bias}(\hat{Q}_3)| \leq \sum_{i=m+1}^{2^{J_*}m} \min\left(\frac{2\tau_i}{n}, \theta_i^2\right) + \sum_{i=2^{J_*}m+1}^{\infty} \theta_i^2.$$

The second term in (55) is easy to bound. Note that for any $j \geq 0$, the $L_p$ ball constraint (9) yields for $p < 2$

$$(56) \qquad \left(\sum_{i=2^jm+1}^{2^{j+1}m} \theta_i^2\right)^{1/2} \leq \left(\sum_{i=2^jm+1}^{2^{j+1}m} |\theta_i|^p\right)^{1/p} \leq M 2^{-js} m^{-s}.$$

Hence for all $\theta \in L_p(\alpha, M)$,

$$\sum_{i=2^{J_*}m+1}^{\infty} \theta_i^2 \leq \sum_{j=J_*}^{\infty} M 2^{-2js} m^{-2s} \leq Cn^{-1/2}$$

for some constant $C > 0$.

It remains to bound the first term in (55). Note that it is straightforward to verify that for all $\theta \in L_p(\alpha, M)$ and all $j \geq 1$,

$$(57) \qquad \sum_{i=2^{j-1}m+1}^{2^j m} |\theta_i|^p \leq M^p 2^{ps} 2^{-jps} m^{-ps}.$$

Hence

$$\sum_{i=m+1}^{2^{J_*}m} \min\left(\frac{2\tau_i}{n}, \theta_i^2\right) = \sum_{j=1}^{J_*} \sum_{i=2^{j-1}m+1}^{2^j m} \min\left(\frac{4j}{n}, \theta_i^2\right)$$

$$= \sum_{j=1}^{J_*} \frac{4j}{n} \sum_{i=2^{j-1}m+1}^{2^j m} \min\left(1, \theta_i^2 \cdot \frac{n}{4j}\right)$$

$$\leq \sum_{j=1}^{J_*} \frac{4j}{n} \sum_{i=2^{j-1}m+1}^{2^j m} \min\left(1, \left\{\theta_i^2 \cdot \frac{n}{4j}\right\}^{p/2}\right),$$

where the last step follows from the facts $\min(1, \theta_i^2 \cdot \frac{n}{4j}) \leq 1$ and $\frac{p}{2} \leq 1$. Hence,

$$\sum_{i=m+1}^{2^{J_*}m} \min\left(\frac{2\tau_i}{n}, \theta_i^2\right) \leq \sum_{j=1}^{J_*} \left(\frac{4j}{n}\right)^{1-p/2} \sum_{i=2^{j-1}m+1}^{2^j m} |\theta_i|^p$$

(58)
$$\leq \left\{ M^p 2^{ps+2-p} \sum_{j=1}^{J_*} j^{1-p/2} 2^{-jps} \right\} \cdot m^{-ps} n^{p/2-1}$$

$$\leq C m^{-ps} n^{p/2-1}$$

for some constant $C > 0$, since $\sum_{j=1}^{\infty} j^{1-p/2} 2^{-jps} < \infty$. Hence, with $m = n^{p/(1+2ps)}$,

(59)
$$\sum_{i=m+1}^{2^{J_*}m} \min\left(\frac{2\tau_i}{n}, \theta_i^2\right) \leq C n^{-(2-p/(1+2ps))/2}.$$

Hence for $p < 2$ and $\alpha \leq \frac{1}{2p}$

(60)
$$\text{Bias}^2(\hat{Q}_3) \leq C n^{-(2-p/(1+2ps))}(1 + o(1)).$$

Equations (53) and (60) together yield

$$E_\theta(\hat{Q}_3 - Q(\theta))^2 \leq \text{Bias}^2(\hat{Q}_3) + \text{Var}(\hat{Q}_3) \leq C n^{-(2-p/(1+2ps))}(1 + o(1)). \quad \square$$

PROOF OF THEOREM 2. The proof of Theorem 2 is analogous to that of Theorem 3 and we only give a brief outline here. Set $m = m_2 = \frac{n}{\log n}$. Then

(61)
$$\text{Var}(\hat{Q}_2) \leq \frac{2m}{n^2} + \frac{4\sum_{i=1}^{m} \theta_i^2}{n} + \frac{6\sum_{i=m+1}^{2^{J_*}m} \theta_i^2}{n} + \sum_{j=1}^{J_*} 2^{j-1} m \cdot \frac{4(2j)^{1/2} + 18}{n^2 e^j}$$

$$\leq \frac{4M}{n}(1 + o(1)).$$

The maximum squared bias of $\hat{Q}_2$ is negligible relative to the minimax risk. This can be shown as follows. Same as in (55) we have

(62)
$$|\text{Bias}(\hat{Q}_2)| \leq \sum_{i=m+1}^{2^{J_*}m} \min\left(\frac{2\tau_i}{n}, \theta_i^2\right) + \sum_{i=2^{J_*}m+1}^{\infty} \theta_i^2.$$

With $m = \frac{n}{\log n}$ and $\alpha p > \frac{1}{2}$, equation (58) yields that for some constant $C > 0$

$$\sum_{i=m+1}^{2^{J_*}m} \min\left(\frac{2\tau_i}{n}, \theta_i^2\right) \leq Cm^{-ps}n^{p/2-1} = C\left(\frac{\log n}{n}\right)^{\alpha p}. \tag{63}$$

For the tail term it follows from (56) that

$$\sum_{i=2^{J_*}m+1}^{\infty} \theta_i^2 \leq \sum_{j=J_*}^{\infty} M 2^{-2js} m^{-2s}$$
$$= M(1 - 2^{-2s})^{-1}(2^{J_*}m)^{-2s} \leq Cn^{-1/2}(\log n)^{-2s} \tag{64}$$

for some constant $C > 0$. Equations (63) and (64) yield $\text{Bias}^2(\hat{Q}_2) = o(\frac{1}{n})$ and Theorem 2 follows. □

PROOF OF THEOREM 4. We shall only consider the case $0 < p < 2$ and $\alpha \leq \frac{1}{2p}$ since the information bound given in (17) can otherwise be applied. The main idea is to inscribe a collection of hyperrectangles inside the parameter space. A prior then mixes over the vertices of the hyperrectangles in this collection and a lower bound for the corresponding Bayes risk and hence minimax risk is given.

Let $\Theta_{k,m}$ be the union of the zero vector $\theta_0 = (0, 0, \dots)$ and the collection of vectors which have exactly $k$ nonzero coordinates equal to $\frac{1}{\sqrt{n}}$ in the first $m$ coordinates and are otherwise equal to zero. We shall write $\Theta_m$ for $\Theta_{[m^{1/2}],m}$. Now suppose that $\hat{Q}$ is an estimator which satisfies

$$E_{\theta_0}(\hat{Q} - Q(\theta_0))^2 \leq c\frac{m}{n^2} \tag{65}$$

for some constant $0 < c < \frac{1}{64}e^{1-e}$. We shall now show that in this case

$$\sup_{\theta \in \Theta_m} E_\theta(\hat{Q} - Q(\theta))^2 \geq \left(\frac{1}{4} - 2e^{(e-1)/2}c^{1/2}\right)\frac{m}{n^2}, \tag{66}$$

and hence for some constant $C > 0$,

$$\inf_{\hat{Q}} \sup_{\theta \in \Theta_m} E_\theta(\hat{Q} - Q(\theta))^2 \geq C\frac{m}{n^2}. \tag{67}$$

Let $\psi_\mu$ be the density of a univariate normal distribution with mean $\mu$ and variance $\frac{1}{n}$. Let $\mathcal{I}(k, m)$ be the class of all subsets of $\{1, \dots, m\}$ of $k$ elements and for $I \in \mathcal{I}(k, m)$ let

$$g_I(y_1, \dots, y_m) = \prod_{j=1}^m \psi_{\mu_j}(y_j),$$

where $\mu_j = \frac{1}{\sqrt{n}}\mathbb{1}(j \in I)$. Finally let

$$g = \frac{1}{\binom{m}{k}} \sum_{I \in \mathcal{I}(k,m)} g_I$$

and $f$ be the joint density of $m$ independent normal random variables each with mean 0 and variance $\frac{1}{n}$. Note that a similar mixture prior was used in [1] to give lower bounds in a nonparametric testing problem. Now note that if

$$E_{g_I}\left(\delta - k\frac{1}{n}\right)^2 \leq C$$

for all $I \in \mathcal{I}(k,m)$, then it follows that

$$E_g\left(\delta - k\frac{1}{n}\right)^2 \leq C.$$

We will now apply the constrained risk inequality of Brown and Low [4]. First we need to calculate a chi-squared distance between $f$ and $g$. This is done as follows. Note that

$$\int \frac{g^2}{f} = \frac{1}{\binom{m}{k}^2} \sum_{I \in \mathcal{I}(k,m)} \sum_{I' \in \mathcal{I}(k,m)} \int \frac{g_I g_{I'}}{f}$$

and simple calculations show that

$$\int \frac{g_I g_{I'}}{f} = e^j,$$

where $j$ is the number of points in the set $I \cap I'$. It follows that

$$\int \frac{g^2}{f} = Ee^J,$$

where $J$ has the hypergeometric distribution

$$P(J = j) = \frac{\binom{k}{j}\binom{m-k}{k-j}}{\binom{m}{k}}.$$

Now note that from [16], page 59,

$$P(J = j) \leq \binom{k}{j}\left(\frac{k}{m}\right)^j\left(1 - \frac{k}{m}\right)^{k-j}\left(1 - \frac{k}{m}\right)^{-k}.$$

Now let $k = [m^{1/2}]$. Then for $m \geq 4$,

$$\left(1 - \frac{k}{m}\right)^{-k} \leq 4$$

and hence

$$P(J=j) \leq 4 \binom{k}{j} \left(\frac{k}{m}\right)^j \left(1-\frac{k}{m}\right)^{k-j}.$$

Consequently

$$\int \frac{g^2}{f} = Ee^J \leq 4\left(1+(e-1)\frac{k}{m}\right)^k \leq 4e^{e-1}.$$

It now follows from the constrained risk inequality in [4] that if

(68) $$E_f(\hat{Q} - Q(\theta_0))^2 \leq c\frac{m}{n^2},$$

then

(69) $$\begin{aligned} E_g\left(\hat{Q} - \frac{k}{n}\right)^2 &\geq \frac{k^2}{n^2} - 4\frac{k}{n}e^{(e-1)/2}c^{1/2}\frac{m^{1/2}}{n} \\ &\geq (1 - 8e^{(e-1)/2}c^{1/2})\frac{k^2}{n^2} \\ &\geq \left(\frac{1}{4} - 2e^{(e-1)/2}c^{1/2}\right)\frac{m}{n^2}. \end{aligned}$$

Hence (67) holds.

It is now easy to check that $\Theta_m$ as defined above is contained in the $L_p$ ball $L_p(\alpha, M)$ when $m = Cn^{p/(1+2ps)}$ for sufficiently small constant $C > 0$. Hence it directly follows from (67) that

(70) $$\inf_{\delta} \sup_{\theta \in L_p(\alpha,M)} E_\theta(\delta - Q(\theta))^2 \geq \inf_{\delta} \sup_{\theta \in \Theta_m} E_\theta(\delta - Q(\theta))^2 \\ \geq Cn^{p/(1+2ps)-2}. \qquad \square$$

**Acknowledgment.** We wish to thank a referee for an extremely detailed and useful report which led to an improvement in some of our earlier results and also helped with the presentation of the paper.

DEPARTMENT OF STATISTICS
THE WHARTON SCHOOL
UNIVERSITY OF PENNSYLVANIA
PHILADELPHIA, PENNSYLVANIA 19104-6340
USA
E-MAIL: tcai@wharton.upenn.edu
       lowm@wharton.upenn.edu